\documentclass[preprint,12pt]{elsarticle}

\usepackage{amssymb}
\usepackage{amsmath}
\usepackage{physics}
\usepackage{mathtools}
\usepackage{bm}
\usepackage{booktabs}
\usepackage{placeins}
\usepackage{float}

\journal{Journal name}

\begin{document}
\begin{frontmatter}

\title{An accelerated direct solver for scalar wave scattering by multiple transmissive inclusions in two dimensions}

\author[label1]{Yasuhiro Matsumoto\corref{cor1}}
\ead{matsumoto@cii.isct.ac.jp}
\cortext[cor1]{Corresponding author}
\affiliation[label1]{organization={Center for Information Infrastructure, Institute of Science Tokyo},
            addressline={2-12-1-I8-21, Ookayama, Meguro-ku},
            city={Tokyo},
            postcode={152-8550},
            country={Japan}}

\begin{abstract}
This paper discusses a fast direct solver using boundary integral equations for Helmholtz transmission problems involving multiple inclusions in two dimensions. Efficiently addressing scattering problems in the presence of numerous inclusions remains a key challenge for various practical applications. For problems involving a large number of scatterers, the number of iterations in Krylov subspace methods is known to increase significantly. This occurs even when using second-kind boundary integral equations, which are typically recognized for their rapid convergence. We consider a fast direct solver as an alternative, an approach that has been less commonly explored for transmission problems with disjoint multiple inclusions. The low-rank approximation based on the proxy method achieve speedup by calculating interactions between disjoint scatterers without the terms derived from the internal integral representation. Notably, this advantage applies to the Poggio--Miller--Chang--Harrington--Wu--Tsai (PMCHWT) formulation but breaks down in the Burton--Miller case. Numerical examples demonstrate that the proposed solver can compress the system of linear algebraic equations to a size of $O(\omega D)$, where $\omega$ is the frequency of the incident wave and $D$ is the diameter of the (smallest) bounding box enclosing the multiple inclusions. The total computational cost scales as $O(N^{1.5})$ $(= O(\sqrt{N}^3))$ at most for a fixed $\omega$ when the inclusions are arranged on a grid. Moreover, the PMCHWT formulation, that omits the interior term in the proxy method, is approximately six times faster than the Burton--Miller formulation when treating each inclusion as a cell. Furthermore, in the same setting, the former can compress the size of the system of linear algebraic equations by half compared to the latter.
\end{abstract}


\end{frontmatter}

\section{Introduction}

Metamaterials are drawing attention in wave scattering problems because they can achieve unique material properties through the precise arrangement and material selection of multiple inclusions (Padilla and Averitt, 2022; Shi et al., 2022).
There is also growing interest in mechanical metamaterials due to their potential to control the scattering of elastic waves as well as scalar waves (Craster et al., 2023).
A computationally efficient numerical approach is essential for solving scattering problems characterized by a large number of inclusions.
Furthermore, boundary integral equation methods (BIEMs) are particularly promising due to their inherent ability to handle radiation conditions in unbounded domains (Nédélec, 2001).
The fast multipole method (FMM) for BIEM continues to be developed to efficiently handle large inclusions, for example by exploiting periodicity (Smith and Abrahams, 2022) and T-matrix formulations (Hawkins et al., 2024).
On the other hand, it is reported that for problems involving a large number of scatterers,
the number of iterations in Krylov subspace methods increases significantly.
This occurs even when using second-kind boundary integral equations (Bremer et al., 2015), which are typically recognized for their rapid convergence.
To further accelerate convergence even in such cases, several strategies have been proposed, including using the FMM itself as a preconditioner (Carpentieri et al., 2005) or employing the LU decomposition of a low-precision $\mathcal{H}^2$-matrix (Spendlhofer et al., 2025), however fast direct solvers (Martinsson, 2019) are also an attractive alternative.
A key advantage of fast direct solvers is that their computational cost remains generally independent of the condition number of the system of linear algebraic equations (Greengard et al., 2009).

While various fast direct solvers have been proposed, this study focuses on solvers based on the proxy method, as they are particularly well-suited for boundary integral equations.
The proxy method was originally proposed to handle the boundary value problem of the two-dimensional Laplace equation, where a low-rank approximation of the coefficient matrix is constructed using a local virtual boundary that surrounds a subset of the boundary (Martinsson and Rokhlin, 2005).
It has also been applied to two- and three-dimensional Helmholtz scattering problems (Martinsson and Rokhlin, 2007; Greengard et al., 2009)
and two-dimensional elastic wave scattering problems (Matsumoto and Maruyama, 2025a).
Fast direct solvers based on the proxy method combined with the quadrature by expansion have also been proposed (Fikl and Klöckner, 2026).
There are many studies applying the proxy methods to scattering problems involving multi-layered media and periodic structures (Gillman, A. and Barnett, 2013; Greengard et al., 2014; Cho and Barnett, 2015; Zhang and Gillman, 2021).
Given that metamaterials are not strictly limited to layered or periodic structures but can also be composed of numerous disjoint scatterers with arbitrary shapes (Matsushima et al., 2022), there is significant room to investigate fast direct solvers for multiple scattering.
Multiple scattering problems are typically addressed within the context of iterative solvers (Gimbutas and Greengard, 2013; Lai and Li, 2019; Hawkins et al., 2024).
However, iterative solvers suffer from the aforementioned issue of slow convergence in multiple scattering problems (Bremer et al., 2015). Therefore, fast direct solvers are expected to offer an advantage even in such scenarios, particularly due to their ability to efficiently handle multiple right-hand sides.
Problems with multiple right-hand sides arise, for example, when varying the position of point sources or the incident angle of a plane wave.
However, fast direct solvers for multiple transmission scattering are rare, whereas their counterparts for exterior or interior problems using the hierarchical off-diagonal low-rank (HODLR) matrix decomposition have been proposed (Bremer et al., 2015).
While fast direct solvers for transmission problems have been developed based on volume integral equations (Gopal and Martinsson, 2022), they lack the dimensional reduction benefits offered by BIEMs.
Consequently, fast direct solvers for multiple transmission scattering problems based on BIEMs should be studied.
In fact, the choice of the appropriate boundary integral equation and the methodology for its low-rank approximation remain non-trivial in this problem.

This work proposes a Poggio--Miller--Chang--Harrington--Wu--Tsai (PMCHWT) formulated fast direct solver for multiple transmission scattering in two dimensions.
We discuss two types of boundary integral equations that directly handle unknown functions. The first type is formulated as a sum of layer potentials with respect to the exterior and interior regions, represented by the PMCHWT formulation (Poggio and Miller, 1973; Chang and Harrington, 1977). The second type involves a formulation where the exterior and interior layer potentials are separated, with the Burton--Miller formulation (Burton and Miller, 1971) chosen as a representative example.
For cases involving a single inclusion, it has been reported that the Burton--Miller formulation is more advantageous in terms of computational time than the PMCHWT formulation, as it can be described using fewer layer potentials (Matsumoto and Maruyama, 2025b).
However, in multiple scattering problems, the calculation of interactions between different scatterers becomes dominant. We point out that while the PMCHWT formulation can leverage this situation to construct a fast direct solver with reduced computational cost, the same construction fails for the Burton--Miller formulation.
In particular, for two-dimensional multiple scattering with smooth boundaries, high-order discretization methods allow us to keep the degrees of freedom per scatterer low.
By identifying each individual scatterer as a cell within the fast direct solver, the acceleration of low-rank approximations via the proxy method is maximized.
Our numerical results indicate that the PMCHWT formulation outperforms the Burton--Miller formulation by a factor of six in terms of computational time when scatterers are distributed on a grid. Additionally, the PMCHWT based approach achieves a 50\% greater reduction in the compressed system's degrees of freedom compared to the solver based on the Burton--Miller formulation.

The remainder of this manuscript is organized as follows. Section 2 describes the transmission problems involving multiple inclusions and the corresponding boundary integral equations. The Nystr\"om discretization for these equations is also presented in this section. Section 3 discusses the proposed fast direct solver, where we point out why the Burton--Miller formulation might be unsuitable for analyzing multiple scatterings. Section 4 demonstrates the performance of the proposed solver through numerical examples, and Section 5 concludes the paper.

\section{Boundary integral formulations for transmission problems with multiple inclusions}
\subsection{Boundary integral equations}
We first state the multiple scattering problems.
Let $\Omega^{-}$ be an union of $m$ disjoint finite domains $\Omega_1, \Omega_2, \ldots \Omega_m$ in $\mathbb{R}^2$.
For $j = 1, 2, \ldots, m$, assume that the boundary $\Gamma_j$ of each $\Omega_j$ is a Jordan closed curve and is smooth.
Let $\Omega^{+} = \mathbb{R}^2 \setminus \overline{\Omega^-}$ be connected.
We assume that there exits an incident wave $u^I$ in $\Omega^{+}$.
Let $k^{\pm} = \omega \sqrt{\varepsilon^{\pm}}$ be the wave number in $\Omega^{\pm}$ for the angular frequency $\omega$ of the incident wave $u^I$
and a material constant $\varepsilon^{\pm}$ in $\Omega^{\pm}$.
Assume that $u^{S \pm}$ is the solution to the multiple transmission scattering problem in $\Omega^{\pm}$,
which satisfies the radiation condition in $\Omega^{+}$.
That is, $u^{S \pm}$ satisfies the Helmholtz equation
\begin{equation}
  \qty{\Delta + (k^{\pm})^2} u^{S \pm}(x) = 0, \quad x \in \Omega^{\pm}.
\end{equation}
The regularity of $u^{S \pm}$ is assumed to be at least $u^{S \pm} \in C^2(\Omega^\pm) \cap C^{1, \gamma}(\overline{\Omega^{\pm}})$ for $0 < \gamma < 1$.
We define the total field $u_j$ on the boundary $\Gamma_j$ by
\begin{equation}
  u_{j}(x) := u^{S+}(x) + u^I(x) = u^{S-}(x), \quad x \in \Gamma_j, \quad j = 1, 2, \ldots, m, \label{eq:pressure}
\end{equation}
and we define its scaled normal derivative $q_j$ on the boundary $\Gamma_j$ by
\begin{equation}
  q_j(x) := \frac{1}{\varepsilon^+} \qty(\pdv{u^{S+}}{n}\qty(x) + \pdv{u^I}{n}\qty(x))
  = \frac{1}{\varepsilon^-} \pdv{u^{S-}}{n}\qty(x), \quad x \in \Gamma_j, \quad j = 1, 2, \ldots, m, \label{eq:flux}
\end{equation}
where $n(z)$ is the unit normal vector at $z \in \bigcup_{j = 1}^{m} \Gamma_j$ pointing toward $\Omega^{+}$, and $\pdv{f}{n}$ represents the normal derivative of a function $f$.
Here, we derive the above two relations by applying the transmission conditions of the multiple scattering problem.
Similarly, let $u_j^{I}$ and $q_j^{I}$ be the incident wave on $\Gamma_j$ and its scaled normal derivative, respectively.

  In the context of acoustic scattering, the material parameter $\varepsilon^\pm$ represents the ratio of density to bulk modulus for the media occupying $\Omega^\pm$.
  The function $u_j$ in the boundary condition \eqref{eq:pressure} corresponds to the sound pressure on the boundary of inclusions.
  This condition \eqref{eq:pressure} means the equilibrium of pressure across the interface.
  Similarly, the function $q_j$ in the boundary condition \eqref{eq:flux} corresponds to the normal particle velocity on the boundary.
  This condition \eqref{eq:flux} ensures continuity of displacement,
implying that no delamination occurs at the interface.

We next formulate two boundary integral equations corresponding to the Helmholtz multiple transmission problem.
The fundamental solution of two-dimensional Helmholtz equation $G^{\pm} (x, y)$ is expressed as
\begin{equation}
  G^{\pm}(x,y) = H_{0}^{(1)}(k^{\pm} |x - y|), \quad x, y \in \mathbb{R}^2, \quad x \neq y,
\end{equation}
where $|\cdot|$ is the Euclidean norm in $\mathbb{R}^2$,
and $H_{0}^{(1)}$ is the zeroth-order Hankel function of the first kind.
Letting $v$ be a function defined on $\bigcup_{j = 1}^{m} \Gamma_j$,
the single layer potential $S_{ij}^{\pm}$,
the double layer potential $D_{ij}^{\pm}$,
and their normal derivatives $D_{ij}^{*\pm}$ and $N_{ij}^{\pm}$
with respect to $\Omega^{\pm}$ are defined by
\begin{equation}
  S_{ij}^{\pm} v(x) = \int_{\Gamma_j} G^{\pm} (x, y) v(y) ds(y), \quad x \in \Gamma_i, \quad j = 1, 2, \ldots, m, \label{eq:slp}
\end{equation}
\begin{equation}
  D_{ij}^{\pm} v(x) = \int_{\Gamma_j} \pdv{G^{\pm} (x, y)}{n(y)} v(y) ds(y), \quad x \in \Gamma_i, \quad j = 1, 2, \ldots, m, \label{eq:dlp}
\end{equation}
\begin{equation}
  D_{ij}^{*\pm} v(x) = \int_{\Gamma_j} \pdv{G^{\pm} (x, y)}{n(x)} v(y) ds(y), \quad x \in \Gamma_i, \quad j = 1, 2, \ldots, m, \label{eq:d_slp}
\end{equation}
\begin{equation}
  N_{ij}^{\pm} v(x) = \mathrm{FP}\int_{\Gamma_j} \pdv{G^{\pm} (x, y)}{n(x)}{n(y)} v(y) ds(y), \quad x \in \Gamma_i, \quad j = 1, 2, \ldots, m, \label{eq:d_dlp}
\end{equation}
respectively.
Here, $\mathrm{FP \int}$ stands for the Hadamard finite part of a divergent integral.
Then, the PMCHWT formulated boundary integral equation is expressed as
\begin{center}
  \makebox[\textwidth][c]{
    \begin{minipage}{1.5\textwidth}
      {\footnotesize
        \begin{equation}
          \begin{aligned}
            \mqty(
            -\qty(\varepsilon^{+} S_{11}^{+} + \varepsilon^{-} S_{11}^{-}) & D_{11}^{+} + D_{11}^{-} & -\qty(\varepsilon^{+} S_{12}^{+} + \varepsilon^{-} S_{12}^{-}) & D_{12}^{+} + D_{12}^{-} & \cdots & -\qty(\varepsilon^{+} S_{1m}^{+} + \varepsilon^{-}S_{1m}^{-}) & D_{1m}^{+} + D_{1m}^{-} \\
            -\qty(D_{11}^{*+} + D_{11}^{*-}) & \frac{1}{\varepsilon^+}N_{11}^{+} + \frac{1}{\varepsilon^-}N_{11}^{-} & -\qty(D_{12}^{*+} + D_{12}^{*-}) & \frac{1}{\varepsilon^+}N_{12}^{+} + \frac{1}{\varepsilon^-}N_{12}^{-} & \cdots & -\qty(D_{1m}^{*+} + D_{1m}^{*-}) & \frac{1}{\varepsilon^+}N_{1m}^{+} + \frac{1}{\varepsilon^-}N_{1m}^{-} \\
            -\qty(\varepsilon^{+} S_{21}^{+} + \varepsilon^{-} S_{21}^{-}) & D_{21}^{+} + D_{21}^{-} & -\qty(\varepsilon^{+} S_{22}^{+} + \varepsilon^{-} S_{22}^{-}) & D_{22}^{+} + D_{22}^{-} & \cdots & -\qty(\varepsilon^{+} S_{2m}^{+} + \varepsilon^{-}S_{2m}^{-}) & D_{2m}^{+} + D_{2m}^{-} \\
            -\qty(D_{21}^{*+} + D_{21}^{*-}) & \frac{1}{\varepsilon^+}N_{21}^{+} + \frac{1}{\varepsilon^-}N_{21}^{-} & -\qty(D_{22}^{*+} + D_{22}^{*-}) & \frac{1}{\varepsilon^+}N_{22}^{+} + \frac{1}{\varepsilon^-}N_{22}^{-} & \cdots & -\qty(D_{2m}^{*+} + D_{2m}^{*-}) & \frac{1}{\varepsilon^+}N_{2m}^{+} + \frac{1}{\varepsilon^-}N_{2m}^{-} \\
            \vdots & \vdots & \vdots & \vdots & \ddots & \vdots & \vdots \\
            -\qty(\varepsilon^{+} S_{m1}^{+} + \varepsilon^{-} S_{m1}^{-}) & D_{m1}^{+} + D_{m1}^{-} & -\qty(\varepsilon^{+} S_{m2}^{+} + \varepsilon^{-} S_{m2}^{-}) & D_{m2}^{+} + D_{m2}^{-} & \cdots & -\qty(\varepsilon^{+} S_{mm}^{+} + \varepsilon^{-}S_{mm}^{-}) & D_{mm}^{+} + D_{mm}^{-} \\
            -\qty(D_{m1}^{*+} + D_{m1}^{*-}) & \frac{1}{\varepsilon^+}N_{m1}^{+} + \frac{1}{\varepsilon^-}N_{m1}^{-} & -\qty(D_{m2}^{*+} + D_{m2}^{*-}) & \frac{1}{\varepsilon^+}N_{m2}^{+} + \frac{1}{\varepsilon^-}N_{m2}^{-} & \cdots & -\qty(D_{mm}^{*+} + D_{mm}^{*-}) & \frac{1}{\varepsilon^+}N_{mm}^{+} + \frac{1}{\varepsilon^-}N_{mm}^{-} 
            ) \\
            \mqty(
            q_{1} \\
            u_{1} \\
            q_{2} \\
            u_{2} \\
            \vdots \\
            q_{m} \\
            u_{m} \\
            )
            =
            \mqty(
            -u_{1}^{I} \\
            -q_{1}^{I} \\
            -u_{2}^{I} \\
            -q_{2}^{I} \\
            \vdots \\
            -u_{m}^{I} \\
            -q_{m}^{I} \\
            ).
          \end{aligned}
          \label{eq:pmchwt_all_bie}
        \end{equation}
      }
    \end{minipage}
  }
\end{center}
By exploiting the fact that layer potentials arising from the interior region 
vanishes in the calculation of the interaction between different inclusions, the above equation can be simplified as
\begin{center}
  \makebox[\textwidth][c]{
    \begin{minipage}{1.4\textwidth}
      {\footnotesize
        \begin{equation}
          \begin{aligned}
            \mqty(
            -\qty(\varepsilon^{+} S_{11}^{+} + \varepsilon^{-} S_{11}^{-}) & D_{11}^{+} + D_{11}^{-} & -\varepsilon^{+} S_{12}^{+}  & D_{12}^{+}  & \cdots & -\varepsilon^{+} S_{1m}^{+} & D_{1m}^{+}  \\
            -\qty(D_{11}^{*+} + D_{11}^{*-}) & \frac{1}{\varepsilon^+}N_{11}^{+} + \frac{1}{\varepsilon^-}N_{11}^{-} & -D_{12}^{*+} & \frac{1}{\varepsilon^+}N_{12}^{+} & \cdots & -D_{1m}^{*+} & \frac{1}{\varepsilon^+}N_{1m}^{+}  \\
            -\varepsilon^{+} S_{21}^{+}  & D_{21}^{+}  & -\qty(\varepsilon^{+} S_{22}^{+} + \varepsilon^{-} S_{22}^{-}) & D_{22}^{+} + D_{22}^{-} & \cdots & -\varepsilon^{+} S_{2m}^{+}  & D_{2m}^{+} \\
            -D_{21}^{*+}  & \frac{1}{\varepsilon^+}N_{21}^{+}  & -\qty(D_{22}^{*+} + D_{22}^{*-}) & \frac{1}{\varepsilon^+}N_{22}^{+} + \frac{1}{\varepsilon^-}N_{22}^{-} & \cdots & -D_{2m}^{*+} & \frac{1}{\varepsilon^+}N_{2m}^{+} \\
            \vdots & \vdots & \vdots & \vdots & \ddots & \vdots & \vdots \\
            -\varepsilon^{+} S_{m1}^{+}  & D_{m1}^{+} & -\varepsilon^{+} S_{m2}^{+}  & D_{m2}^{+} & \cdots & -\qty(\varepsilon^{+} S_{mm}^{+} + \varepsilon^{-}S_{mm}^{-}) & D_{mm}^{+} + D_{mm}^{-} \\
            -D_{m1}^{*+} & \frac{1}{\varepsilon^+}N_{m1}^{+}  & -D_{m2}^{*+} & \frac{1}{\varepsilon^+}N_{m2}^{+} & \cdots & -\qty(D_{mm}^{*+} + D_{mm}^{*-}) & \frac{1}{\varepsilon^+}N_{mm}^{+} + \frac{1}{\varepsilon^-}N_{mm}^{-} 
            ) \\
            \mqty(
            q_{1} \\
            u_{1} \\
            q_{2} \\
            u_{2} \\
            \vdots \\
            q_{m} \\
            u_{m} \\
            )
            =
            \mqty(
            -u_{1}^{I} \\
            -q_{1}^{I} \\
            -u_{2}^{I} \\
            -q_{2}^{I} \\
            \vdots \\
            -u_{m}^{I} \\
            -q_{m}^{I} \\
            ).
          \end{aligned}
          \label{eq:pmchwt_omit_bie}
        \end{equation}
      }
    \end{minipage}
  }
\end{center}
We refer to the formulations of \eqref{eq:pmchwt_all_bie} and \eqref{eq:pmchwt_omit_bie} as ``PMCHWT (All)'' and ``PMCHWT Omit'' in this work, respectively.
The Burton--Miller formulated boundary integral equation,
which provides an alternative formulation of the same problem, is given by
\begin{center}
  \makebox[\textwidth][c]{
    \begin{minipage}{1.4\textwidth}
      {\footnotesize
        \begin{equation}
          \begin{aligned}
            \mqty(
            \begin{aligned} \qty(D_{11}^{+} - 1/2) \\ + \alpha N_{11}^{+} \end{aligned} & \begin{aligned} &-\varepsilon^{+} S_{11}^{+} \\ &- \alpha \varepsilon^{+} \qty(D_{11}^{*+} + 1/2) \end{aligned} & D_{12}^{+} + \alpha N_{12}^{+} & -\varepsilon^{+}\qty(S_{12}^{+} + \alpha D_{12}^{*+}) & \cdots & D_{1m}^{+} + \alpha N_{1m}^{+} & -\varepsilon^{+}\qty(S_{1m}^{+} + \alpha D_{1m}^{*+}) \\
            D_{11}^{-} + 1/2 & -\varepsilon^- S_{11}^{-} & D_{12}^{-} & -\varepsilon^- S_{12}^{-} & \cdots & D_{1m}^{-} & -\varepsilon^- S_{1m}^{-} \\
            D_{21}^{+} + \alpha N_{21}^{+} & -\varepsilon^{+}\qty(S_{21}^{+} + \alpha D_{21}^{*+}) & \begin{aligned} \qty(D_{22}^{+} - 1/2) \\ + \alpha N_{22}^{+} \end{aligned} & \begin{aligned} &-\varepsilon^{+} S_{22}^{+} \\ &- \alpha \varepsilon^{+} \qty(D_{22}^{*+} + 1/2) \end{aligned} & \cdots & D_{2m}^{+} + \alpha N_{2m}^{+} & -\varepsilon^{+}\qty(S_{2m}^{+} + \alpha D_{2m}^{*+}) \\
            D_{21}^{-} & -\varepsilon^- S_{21}^{-} & D_{22}^{-} + 1/2 & -\varepsilon^- S_{22}^{-} & \cdots & D_{2m}^{-} & -\varepsilon^- S_{2m}^{-} \\
            \vdots & \vdots & \vdots & \vdots & \ddots & \vdots & \vdots \\
            D_{m1}^{+} + \alpha N_{m1}^{-} & -\varepsilon^{+}\qty(S_{m1}^{+} + \alpha D_{m1}^{*+}) & D_{m2}^{+} + \alpha N_{m2}^{-} & -\varepsilon^{+}\qty(S_{m2}^{+} + \alpha D_{m2}^{*+}) & \cdots & \begin{aligned} \qty(D_{mm}^{+} - 1/2) \\ + \alpha N_{mm}^{-} \end{aligned} & \begin{aligned} &-\varepsilon^{+} S_{mm}^{+} \\ &- \alpha \varepsilon^{+} \qty(D_{mm}^{*+} + 1/2) \end{aligned} \\
            D_{m1}^{-} & -\varepsilon^- S_{m1}^{-} & D_{m2}^{-} & -\varepsilon^- S_{m2}^{-} & \cdots & D_{mm}^{-} + 1/2 & -\varepsilon^- S_{mm}^{-}
            ) \\
            \mqty(
            u_{1} \\
            q_{1} \\
            u_{2} \\
            q_{2} \\
            \vdots \\
            u_{m} \\
            q_{m} \\
            )
            =
            \mqty(
            -\qty(u_{1}^{I} + \alpha q_{1}^{I}) \\
            0 \\
            -\qty(u_{2}^{I} + \alpha q_{2}^{I}) \\
            0 \\
            \vdots \\
            -\qty(u_{m}^{I} + \alpha q_{m}^{I}) \\
            0
            ),
          \end{aligned}
          \label{eq:bm_bie}
        \end{equation}
      }
    \end{minipage}
  }
\end{center}
where $\alpha = i/k^+$ is a constant for Burton--Miller method.
Here $i$ is the imaginary unit.
We refer to the formulation of \eqref{eq:bm_bie} as ``BM'' in this work.
Although a simplified Burton--Miller formulation could be considered, similar to the PMCHWT formulation, we employ this naive form to construct a fast direct solver for reasons discussed later.

\subsection{Discretization}
We briefly denote the discretization for the boundary integral equations \eqref{eq:pmchwt_all_bie}, \eqref{eq:pmchwt_omit_bie} and \eqref{eq:bm_bie}.
We use a Nystr\"om discretization method with the zeta corrected quadrature (Wu and Martinsson, 2021 and 2023), applying $n$ quadrature points for each boundary $\Gamma_j$, $j = 1, 2, \ldots, m$.
The zeta corrected quadrature is the high-order local-corrected quadrature for the layer potentials \eqref{eq:slp}, \eqref{eq:dlp}, \eqref{eq:d_slp} and \eqref{eq:d_dlp}.
This study employs 41 local correction points, resulting in a quadrature scheme with an $O(n^{-41})$ convergence rate.
The assumption of smooth inclusion boundaries is necessary to allow the use of a larger number of local correction points.
Let $\bm{S}_{ij}^{\pm}$, $\bm{D}_{ij}^{\pm}$, $\bm{D}_{ij}^{*\pm}$ and $\bm{N}_{ij}^{\pm}$ be the $n \times n$ matrix obtained by discretizing the layer potentials \eqref{eq:slp}, \eqref{eq:dlp}, \eqref{eq:d_slp} and \eqref{eq:d_dlp} representing the interaction between $\Gamma_i$ and $\Gamma_j$, respectively.
Let $A_{ij}$, $\varphi_i$ and $f_i$ be the matrix stemming from the interaction between $\Gamma_i$ and $\Gamma_j$, the vector of the solution corresponding to the quadrature points on $\Gamma_i$, and the vector of the right-hand side corresponding to the quadrature points on $\Gamma_i$, respectively.
Then, boundary integral equations are discretized into a system of linear algebraic equations expressed as
\begin{equation}
  \mqty(
  A_{11} & A_{12} & \cdots & A_{1m} \\
  A_{21} & A_{22} & \cdots & A_{2m} \\
  \vdots & \vdots &\ddots &\vdots \\
  A_{m1} & A_{m2} & \cdots & A_{mm}
  )
  \mqty(
  \varphi_1 \\
  \varphi_2 \\
  \vdots \\
  \varphi_m
  )
  =
  \mqty(
  f_1 \\
  f_2 \\
  \vdots \\
  f_m
  ).
  \label{eq:block_linear}
\end{equation}
In the PMCHWT (All) formulation \eqref{eq:pmchwt_all_bie},
$A_{ij}$, $\varphi_i$ and $f_i$ are expressed as
\begin{equation}
  A_{ij} = 
  \mqty(
  -\qty(\varepsilon^{+} \bm{S}_{ij}^{+} + \varepsilon^{-} \bm{S}_{ij}^{-}) & \bm{D}_{ij}^{+} + \bm{D}_{ij}^{-} \\
  -\qty(\bm{D}_{ij}^{*+} + \bm{D}_{ij}^{*-}) & \frac{1}{\varepsilon^+} \bm{N}_{ij}^{+} + \frac{1}{\varepsilon^-} \bm{N}_{ij}^{-}
  ),
\end{equation}
\begin{equation}
  \varphi_{i} = \mqty(
  \bm{q}_i \\
  \bm{u}_i \\
  ),
\end{equation}
\begin{equation}
  f_{i} = \mqty(
  -\bm{u}_i^I \\
  -\bm{q}_i^I \\
  ),
\end{equation}
for $i, j = 1, 2, \ldots, m$,
where $\bm{u}_i$, $\bm{q}_i$, $\bm{u}_i^I$ and $\bm{q}_i^I$ are discretized $u_i$, $q_i$, $u_i^I$ and $q_i^I$.
In the PMCHWT (Omit) formulation \eqref{eq:pmchwt_omit_bie},
$A_{ij}$, $\varphi_i$ and $f_i$ are expressed as
\begin{equation}
  A_{ij} = 
  \begin{dcases}
    \mqty(
    -\qty(\varepsilon^{+} \bm{S}_{ij}^{+} + \varepsilon^{-} \bm{S}_{ij}^{-}) & \bm{D}_{ij}^{+} + \bm{D}_{ij}^{-} \\
    -\qty(\bm{D}_{ij}^{*+} + \bm{D}_{ij}^{*-}) & \frac{1}{\varepsilon^+} \bm{N}_{ij}^{+} + \frac{1}{\varepsilon^-} \bm{N}_{ij}^{-}
    ) &(i = j), \\
    \mqty(
    -\varepsilon^{+} \bm{S}_{ij}^{+} & \bm{D}_{ij}^{+}  \\
    -\bm{D}_{ij}^{*+} & \frac{1}{\varepsilon^+} \bm{N}_{ij}^{+}
    ) &(i \neq j), \\
  \end{dcases}
\end{equation}
\begin{equation}
  \varphi_{i} = \mqty(
  \bm{q}_i \\
  \bm{u}_i \\
  ),
\end{equation}
\begin{equation}
  f_{i} = \mqty(
  -\bm{u}_i^I \\
  -\bm{q}_i^I \\
  ),
\end{equation}
for $i, j = 1, 2, \ldots, m$.
In the BM formulation \eqref{eq:bm_bie},
$A_{ij}$, $\varphi_i$ and $f_i$ are expressed as
\begin{equation}
  A_{ij} = 
  \begin{dcases}
    \mqty(
    \bm{D}_{ij}^{+} - I/2 + \alpha \bm{N}_{ij}^{+} & -\varepsilon^{+} \bm{S}_{ij}^{+} - \alpha \varepsilon^{+} \qty(\bm{D}_{ij}^{*+} + I/2) \\
    \bm{D}_{ij}^{-} + I/2 & -\varepsilon^- \bm{S}_{ij}^{-}
    ) &(i = j), \\
    \mqty(
    \bm{D}_{ij}^{+} + \alpha \bm{N}_{ij}^{+} & -\varepsilon^{+} \bm{S}_{ij}^{+} - \alpha \varepsilon^{+} \bm{D}_{ij}^{*+} \\
    \bm{D}_{ij}^{-}  & -\varepsilon^- \bm{S}_{ij}^{-}
    ) &(i \neq j), \\
  \end{dcases}
\end{equation}
\begin{equation}
  \varphi_{i} = \mqty(
  \bm{u}_i \\
  \bm{q}_i
  ),
\end{equation}
\begin{equation}
  f_{i} = \mqty(
  -\bm{u}_i^I -\alpha \bm{q}_i^I \\
  0
  ),
\end{equation}
for $i, j = 1, 2, \ldots, m$, where $I$ is the $n \times n$ identity matrix.

\section{Fast direct solver}
Since the fast direct solver used in this work is fundamentally based on our previous studies (Matsumoto and Matsushima, 2025),
we will focus our discussion on the specific differences introduced here for multiple scattering problems.
A tree structure is employed to hierarchically partition the boundary set $\bigcup_{j = 1}^{m} \Gamma_j$.
The fast direct solver primarily consists of two phases: a compression step performed upward through the tree structure and a conversion step executed downward.
The point to be noted lies in the former.

This work assumes that the coefficient matrix of \eqref{eq:block_linear} is weakly admissible,
which implies that every off-diagonal block possesses a low-rank structure.
This assumption is justified in two-dimensional problems (Ma et al., 2022).
In the following formulation, the rank for low-rank approximations is fixed at $k$ to simplify the presentation.
We assume $k \ll n$, where $n$ is the number of quadrature points of each $\Gamma_j$ ($j = 1, 2, \ldots, m$).
It is further assumed that the geometry is partitioned into cells based on a perfect binary tree structure, where each $\Gamma_j$ itself serves as a cell at the leaf level.
Consequently, at the leaf level, the off-diagonal blocks of the system of linear algebraic equations \eqref{eq:block_linear} can be approximated by
\begin{equation}
  \mqty(
  A_{11} & L_{1}S_{12}R_{2} & \cdots & L_{1}S_{1m}R_{m} \\
  L_{2}S_{21}R_{1} & A_{22} & \cdots & L_{2}S_{2m}R_{m} \\
  \vdots & \vdots &\ddots &\vdots \\
  L_{m}S_{m1}R_{1} & L_{m}S_{m2}R_{2} & \cdots & A_{mm}
  )
  \mqty(
  \varphi_1 \\
  \varphi_2 \\
  \vdots \\
  \varphi_m
  )
  =
  \mqty(
  f_1 \\
  f_2 \\
  \vdots \\
  f_m
  ),
  \label{eq:low_rank_block_linear}
\end{equation}
where $S_{ij} \in \mathbb{C}^{2k \times 2k}$ is the so-called skeleton of $A_{ij}$ $(i \neq j)$, while $L_i \in \mathbb{C}^{2n \times 2k}$ and $R_{j} \in \mathbb{C}^{2k \times 2n}$ are interpolation matrices for $S_{ij}$.
Their sub-matrices structures are given by
\begin{equation}
  L_i S_{ij} R_{j} =
  \mqty(
  L_i^1 & \\
  & L_i^2
  )
  \mqty(
  S_{ij}^1 & S_{ij}^2 \\
  S_{ij}^3 & S_{ij}^4
  )
  \mqty(
  R_i^1 & \\
  & R_i^2
  ). \quad i \neq j.
\end{equation}
These $S_{ij}$, $L_i$ and $R_i$ can be computed by the proxy method (Martinsson and Rokhlin, 2005)
and the interpolative decomposition (Cheng et al., 2005).
The proxy method approximates far-field interactions by means of a local virtual boundary.
Since focusing solely on the far-field is inadequate for constructing a weakly admissible low-rank approximation, the proxy method numerically decomposes not only far-field but also near-field interactions with respect to adjacent cells via a virtual boundary.
The skeleton $S_{ij}$ is obtained as a sub-matrix of $A_{ij}$, while $L_i$ and $R_i$ serve as interpolation matrices that reconstruct $A_{ij}$ from $S_{ij}$.
We define $\tilde{A}_i$ by $(R_i A_{ii}^{-1} L_i)^{-1}$.
Then, the system \eqref{eq:low_rank_block_linear} of size $N \times N = 2mn \times 2mn$ can be reduced to the following compressed form as
\begin{equation}
  \mqty(
  \tilde{A}_{1} & S_{12} & \cdots & S_{1m} \\
  S_{21} & \tilde{A}_{2} & \cdots & S_{2m} \\
  \vdots & \vdots &\ddots &\vdots \\
  S_{m1} & S_{m2} & \cdots & \tilde{A}_{m}
  )
  \mqty(
  \psi_1 \\
  \psi_2 \\
  \vdots \\
  \psi_m
  )
  =
  \mqty(
  \tilde{A}_{1}R_{1}A_{11}^{-1}f_1 \\
  \tilde{A}_{2}R_{2}A_{22}^{-1}f_2 \\
  \vdots \\
  \tilde{A}_{m}R_{m}A_{mm}^{-1}f_m
  ),
  \label{eq:compressed}
\end{equation}
where the vectors $\{\psi_i \}_{i = 1}^{m}$ are defined as $\psi_i = R_i \varphi_i$ and introduced as new unknowns.
Because this reduced system has a size of $2mk \times 2mk$, it can be solved efficiently than the original system \eqref{eq:block_linear}.

In the preceding explanation, while the existence of $\tilde{A}_i$ was assumed, it should be noted that when employing the Burton--Miller formulation with the simplified form, the term $(R_i A_{ii}^{-1} L_i)$ becomes singular.
This occurs because this simplified form neglects the layer potential from the interior region in the off-diagonal blocks $A_{ij}$ such as
\begin{equation}
  A_{ij} = 
    \mqty(
    \bm{D}_{ij}^{+} + \alpha \bm{N}_{ij}^{+} & -\varepsilon^{+} \bm{S}_{ij}^{+} - \alpha \varepsilon^{+} \bm{D}_{ij}^{*+} \\
    0 & 0
    ), \quad i \neq j,
\end{equation}
for $i, j = 1, 2, \ldots, m$.
In the simplified Burton--Miller formulation, since the lower blocks of $A_{ij}$ are zero, the corresponding skeleton $S_{ij}$ also contains zero lower blocks. This implies that the interpolation matrix $L_i^2$ always vanishes, leading to the conclusion that $(R_i A_{ii}^{-1} L_i)$ is not invertible.
We therefore have no choice but to use \eqref{eq:bm_bie} to construct the Burton--Miller formulation based fast direct solver.
In contrast, this issue does not arise when handling the PMCHWT (Omit) formulation \eqref{eq:pmchwt_omit_bie}.

After solving \eqref{eq:compressed} for $\{\psi_i \}_{i = 1}^{m}$,
the original solution vectors $\{\varphi_i \}_{i = 1}^{m}$ are recovered by converting $\{ \psi_i \}_{i = 1}^{m}$ as follows:
\begin{equation}
  \mqty(
  \varphi_1 \\
  \varphi_2 \\
  \vdots \\
  \varphi_m
  )
  =
  \mqty(
  A_{11}^{-1}f_1 \\
  A_{22}^{-1}f_2 \\
  \vdots \\
  A_{mm}^{-1}f_m
  )
  -
  \mqty(
  A_{11}^{-1}L_{1} \tilde{A}_{1}R_{1}A_{11}^{-1}f_1 \\
  A_{22}^{-1}L_{2}\tilde{A}_{2}R_{2}A_{22}^{-1}f_2 \\
  \vdots \\
  A_{mm}^{-1}L_{m}\tilde{A}_{m}R_{m}A_{mm}^{-1}f_m
  )
  +
  \mqty(
  A_{11}^{-1}L_{1} \tilde{A}_{1} \psi_1 \\
  A_{22}^{-1}L_{2} \tilde{A}_{2} \psi_2 \\
  \vdots \\
  A_{mm}^{-1}L_{m} \tilde{A}_{m} \psi_m
  ).
\end{equation}
This results in the single-level solver.
At the parent level, the skeleton matrices are formed by merging those from the children, offering the potential for recursive low-rank approximations. By exploiting this property, we can construct a multi-level solver (Martinsson and Rokhlin, 2005).
We use this multi-level algorithm in numerical examples.

Although we used a fixed number of skeletons $k$ in the aforementioned formulation, this rank can be adaptively chosen, and indeed, our implementation supports adaptive ranks.
The adaptive rank is determined during the column-pivoted QR decomposition used for computing the interpolative decomposition; it is decided by checking whether the diagonal components of the resulting upper triangular matrix exceed a predefined threshold $\epsilon$.
However, this adaptive rank approach does not strictly guarantee a low-rank approximation with the accuracy of $\epsilon$, but rather only indirectly controls the precision.
While it is possible to guarantee the decomposition accuracy by incrementally increasing the rank and monitoring the singular values (Martinsson, 2019), we did not adopt this method due to its prohibitive computational cost.

\section{Numerical examples}
\begin{figure}[!tb]
    \centering
    \includegraphics[width=0.5\linewidth]{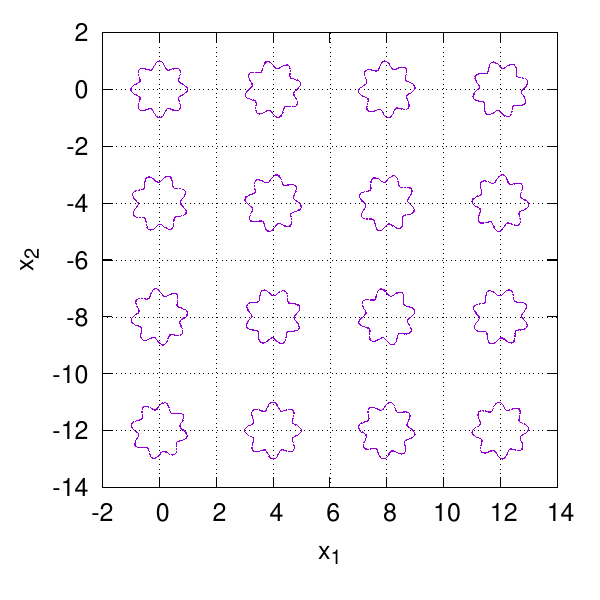}
    \caption{Multiple star-shaped scatterers. This figure corresponds to $4 \times 4 = 16$ scatterers case.}
    \label{fig:star}
\end{figure}

This work assumes some restrictions for inclusions for the numerical experiments.
The inclusions are subject to specific positional constraints: the total number is chosen as a power of 4, and they are aligned on a regular grid.
Figure \ref{fig:star} illustrates an example for the arrangement of inclusions for the case of $4^2 = 16$ inclusions.
In this figure, star-shaped inclusions with a smooth boundary are arranged on a grid while being rotated.
In the numerical experiments, additional cases involving $8^2 = 64$, $16^2 = 256$, $32^2 = 1024$, and $64^2 = 4096$ inclusions of the same shape are also investigated.
Under these restrictions, the diameter $D$ of the bounding box covering the bounded domain is proportional to $\sqrt{N}$,
where $N=nm$ denotes the total number of quadrature points.
Here, $n$ is the number of quadrature points per inclusion and $m$ is the total number of inclusions.
Assuming the compressed system size scales with $O(\omega D)$, consistent with the multipole expansion requirements in the FMM (Rokhlin, 1990), the computational complexity of the proposed fast direct solver for a fixed frequency is expected to be $O(\sqrt{N}^3) = O(N^{1.5})$.
Through several numerical examples, we investigate whether the computational time and the degrees of freedom of the compressed system of linear algebraic equations \eqref{eq:compressed} follow this prediction.

\begin{table}[h]
    \centering
    \caption{Relationship between the number of the star-shaped inclusions and the degrees of freedom ($n = 200$).}
    \label{tb:dof}
    \begin{tabular}{rr}
    \toprule
      Number of inclusions & Degrees of freedom \\
    \midrule
      $16$    & 6400 \\
      $64$    & 25600 \\
      $256$  & 102400 \\
      $1024$ & 409600 \\
      $4096$ & 1638400 \\
    \bottomrule
    \end{tabular}
\end{table}
We describe the general setting for numerical experiments.
We test the performance of the proposed fast direct solver based on PMCHWT (Omit) \eqref{eq:pmchwt_omit_bie}, PMCHWT (all) \eqref{eq:pmchwt_all_bie}, and BM formulations \eqref{eq:bm_bie}.
For the near-field interaction in the proxy method, we consider the inclusions immediately adjacent to the target (above, below, left, and right). This configuration is adopted based on its empirical effectiveness in our experience (Matsumoto, 2025).
The number of quadrature points per inclusion is set to $n = 200$.
The relationship between the number of inclusions and the degrees of freedom in the discrete system is shown in Table \ref{tb:dof}.
A plane wave $u^I = \exp(i k^{+} x_1)$ traveling in the $x_1$ direction is assumed as the incident wave.
For comparison, a conventional BIEM without acceleration is employed.
The same quadrature rule as in the fast direct solver is applied, and the PMCHWT (Omit) and BM formulations are used in the conventional BIEM.
In both the conventional BIEM and the compressed system \eqref{eq:compressed} of the fast direct solver, the resulting linear algebraic equations are solved using a partial pivoted LU decomposition.
A single node of the supercomputer Camphor3 at Kyoto University, which has two 56-core Intel Xeon Max 9480 CPUs (total number of cores: 112; maximum frequency: 1900 MHz) and 128 GB of 3.2 TB/s HBM2e memory, was used for all numerical examples.
In this paper, ``computational time'' refers to the total elapsed time from mesh generation to the acquisition of the numerical solution.
  We note that the mesh generation time accounts for about $0.01 \%$ of the total solver time.
The numerical code was developed in C++ and compiled using the g++ compiler, with parallelization implemented via OpenMP.

\subsection{Verification of conventional boundary integral equation methods} \label{sec:verify}
\begin{figure}[!tb]
    \centering
    \includegraphics[width=1\linewidth]{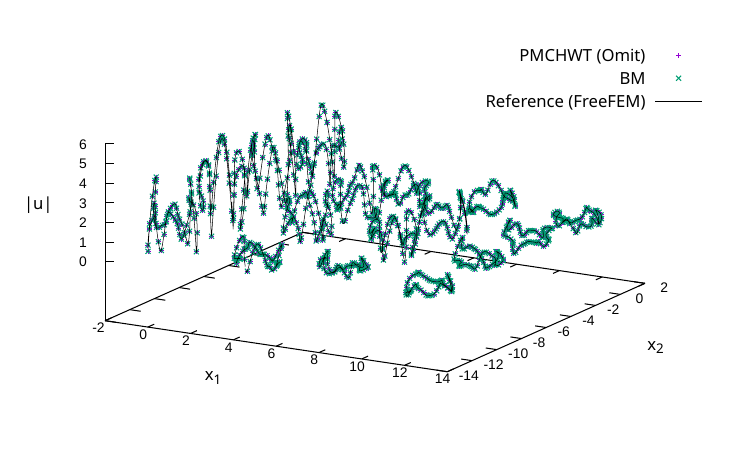}
    \caption{Comparison of the absolute value of the numerical solution $|u(x)| = \sqrt{(\Re(u_j (x)))^2 + (\Im(u_j (x)))^2}$ at quadrature points on the star-shaped boundary $\Gamma_j$ for $j = 1, 2, \ldots, 16$.
        The solid line stands for the reference solution obtained using the BEM module of FreeFEM.
        The markers $+$ and $\times$ correspond to the solutions obtained by the proposed fast direct solver using PMCHWT (Omit) and BM formulations, respectively.
        Both solutions are in good agreement with the reference solution.
    }
    \label{fig:multi}
\end{figure}
\begin{figure}[!tb]
  \centering
  \includegraphics[width=0.8\linewidth]{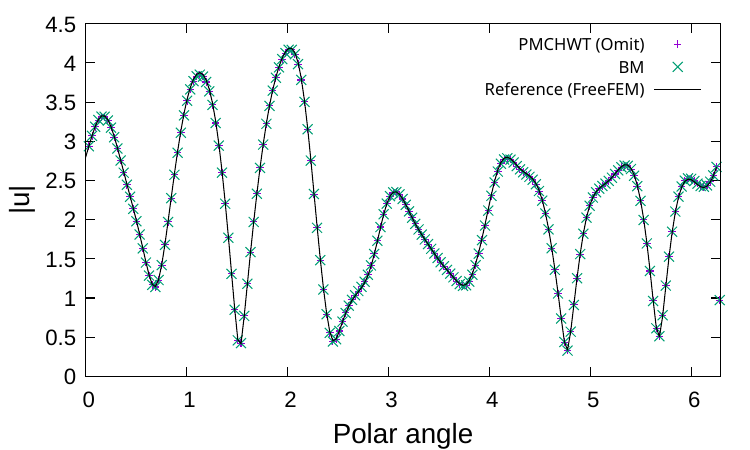}
  \caption{Comparison of the numerical solution on the boundary $\Gamma_1$ of the inclusion centered at $(x_1, x_2) = (0, 0)$.
    The magnitude $|u (x(\theta))| = \sqrt{\qty{\Re(u_1 (x(\theta)))}^2 + \qty{\Im(u_1 (x(\theta)))}^2}$ is plotted at quadrature points, where $\theta$ is the polar angle and $x(\theta)$ represents the parametric boundary.
      The solid line stands for the reference solution obtained using the BEM module of FreeFEM.
      The markers $+$ and $\times$ correspond to the solutions obtained by the proposed fast direct solver using PMCHWT (Omit) and BM formulations, respectively.
      Both solutions are in good agreement with the reference solution.
  }
  \label{fig:multi_u}
\end{figure}
\begin{figure}[!tb]
  \centering
  \includegraphics[width=0.8\linewidth]{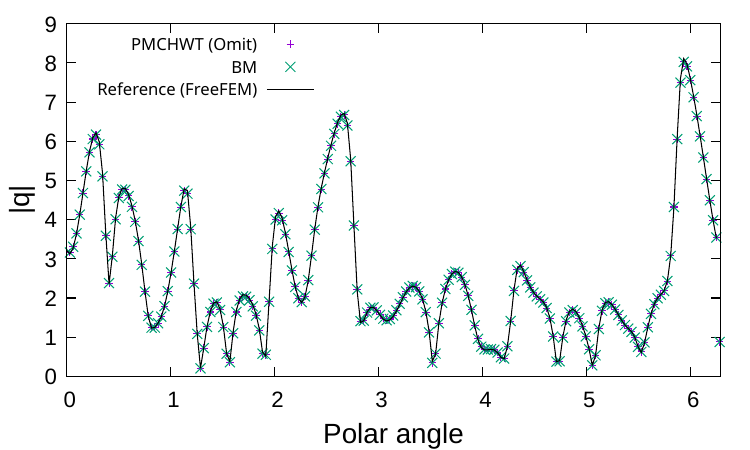}
  \caption{Comparison of the numerical solution on the boundary $\Gamma_1$ of the inclusion centered at $(x_1, x_2) = (0, 0)$.
      The magnitude $|q (x(\theta))| = \sqrt{\qty{\Re(q_1 (x(\theta)))}^2 + \qty{\Im(q_1 (x(\theta)))}^2}$ is plotted at quadrature points, where $\theta$ is the polar angle and $x(\theta)$ represents the parametric boundary.
      The solid line stands for the reference solution obtained using the BEM module of FreeFEM.
      The markers $+$ and $\times$ correspond to the solutions obtained by the proposed fast direct solver using PMCHWT (Omit) and BM formulations, respectively.
      Both solutions are in good agreement with the reference solution.
  }
  \label{fig:multi_q}
\end{figure}
We first verify the implementation of the conventional BIEMs based on both PMCHWT (Omit) \eqref{eq:pmchwt_omit_bie} and BM \eqref{eq:bm_bie} formulations.
As a reference, we employ the numerical solution obtained from the BEM module of FreeFEM (Hecht, 2012) using the BM formulation \eqref{eq:bm_bie}.
In BEM module of FreeFEM, the Galerkin method is employed with P1 test and basis functions.
The calculation conditions used for this verification are shown in Table \ref{tb:verify_condition}.
\begin{table}[H]
    \centering
    \caption{Calculation conditions used in Section \ref{sec:verify}.}
    \label{tb:verify_condition}
    \begin{tabular}{lr}
    \toprule
      Parameter & Value \\
    \midrule
      Angular frequency ($\omega$)    & 3 \\
      Material constant of exterior domain ($\varepsilon^{+}$) & 1 \\
      Material constant of interior domain ($\varepsilon^{-}$) & 4 \\
      Number of inclusions & 16 \\
      Number of quadrature points per inclusion (Ours) & 200 \\
      Number of nodes per inclusion (FreeFEM) & 200 \\
    \bottomrule
    \end{tabular}
\end{table}

Figure \ref{fig:multi} shows the absolute values of the numerical solutions obtained by the proposed conventional BIEMs and FreeFEM. Furthermore, Figures \ref{fig:multi_u} and \ref{fig:multi_q} provide a closer look at the absolute values of the numerical solution and its scaled normal derivative, respectively, for the inclusion centered at $(x_1, x_2) = (0, 0)$.
These figures verify our implementation of the conventional BIEM based on both the PMCHWT (Omit) \eqref{eq:pmchwt_omit_bie} and BM \eqref{eq:bm_bie} formulations for multiple scattering problems.
In the following numerical experiments, this conventional BIEM is used as a reference for accuracy.

\subsection{Increasing the number of inclusions} \label{sec:nlev}
We investigate the performance of the proposed solver when the number of scatterers is increased.
Let $\varphi$ be the numerical solution vector obtained by a specific method.
Its relative 2-norm error, measured against the reference solution $\varphi_{\mathrm{BM}}^{\mathrm{Conv}}$ obtained via the conventional BIEM based on BM formulation, is defined as 
\begin{equation}
  \frac{\norm{\varphi - \varphi_{\mathrm{BM}}^{\mathrm{Conv}}}_2}{\norm{\varphi_{\mathrm{BM}}^{\mathrm{Conv}}}_2}, \label{eq:relative_norm}
\end{equation}
where $\norm{\cdot}_2$ is the 2-norm in $\mathbb{C}^{2N}$.
The calculation conditions used in this subsection are shown in Table \ref{tb:nlev_condition}.
For each formulation, we evaluated the performance of the fast direct solver by varying not only the total number of inclusions but also the threshold for the column-pivoted QR decomposition, which adaptively determines the number of skeleton within the fast direct solver.
Due to computer memory capacity limitations, numerical experiments with 4096 inclusions were conducted only for the PMCHWT (Omit) based fast direct solver with thresholds of $\epsilon = 10^{-9}$ and $10^{-6}$.
\begin{table}[H]
    \centering
    \caption{Calculation conditions used in Section \ref{sec:nlev}.}
    \label{tb:nlev_condition}
    \begin{tabular}{lr}
      \toprule
      Parameter & Value \\
      \midrule
      Angular frequency ($\omega$)    & 5 \\
      Material constant of exterior domain ($\varepsilon^{+}$) & 1 \\
      Material constant of interior domain ($\varepsilon^{-}$) & 5 \\
      Number of inclusions & \{16, 64, 256, 1024, 4096\} \\
      Number of quadrature points per inclusion & 200 \\
      Threshold of the column-pivoted QR decomposition ($\epsilon$) & \{$10^{-12}, 10^{-9}, 10^{-6}$\} \\
      \bottomrule
    \end{tabular}
\end{table}

\begin{figure}[!tb]
  \centering
  \includegraphics[width=0.95\linewidth]{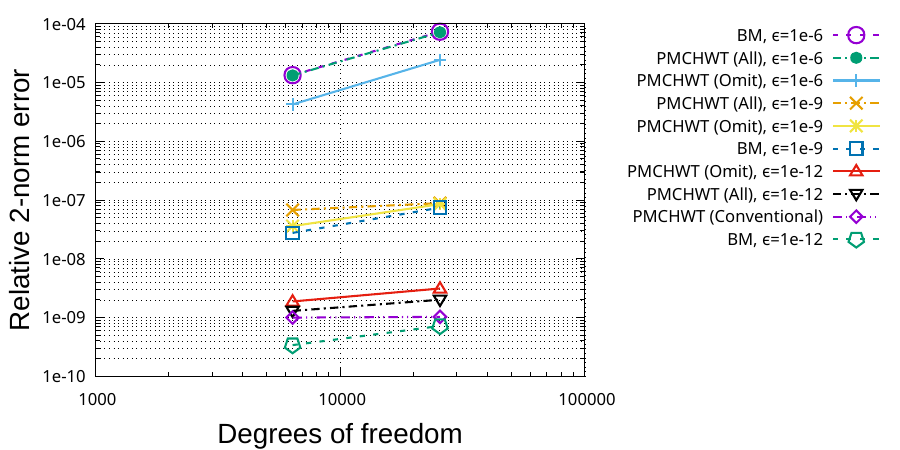}
  \caption{
    Relative 2-norm error defined in \eqref{eq:relative_norm}. Plots are sorted in descending order of their relative errors.
    The labels ``PMCHWT (Omit)'',  ``PMCHWT (All)'' and ``BM'' correspond to the fast direct solver based on \eqref{eq:pmchwt_omit_bie}, \eqref{eq:pmchwt_all_bie} and \eqref{eq:bm_bie}, respectively.
    The labels ``PMCHWT (Conventional)'' correspond to \eqref{eq:pmchwt_omit_bie}, and is solved by a standard partial pivoted LU decomposition.
    In the labels, $\epsilon$ represents the threshold of the column-pivoted QR decomposition in the fast direct solver.
      The results of the PMCHWT (Omit) formulation are represented by solid lines, while those of the other formulations are shown as dashed lines.
      The result of Burton--Miller formulated conventional BIEM is used as the reference.
      In this figure, the angular frequency $\omega$ and the material constants $\varepsilon^\pm$ are fixed.
      It can be observed that PMCHWT (Omit) formulation has negligible impact on accuracy.
  }
  \label{fig:error_nlev}
\end{figure}

Figure \ref{fig:error_nlev} shows the relative 2-norm error \eqref{eq:relative_norm} of each solver.
The solution obtained from the Burton--Miller formulation solved by a standard partial pivoted LU decomposition is employed as the reference.
This figure indicates that the accuracy of the fast direct solver can be controlled through the threshold $\epsilon$.
When $\epsilon = 10^{-12}$ is employed, the fast direct solvers achieve an accuracy on par with the conventional BIEM.
The PMCHWT (All) and PMCHWT (Omit) formulated fast direct solvers yield similar levels of precision.
From Figure \ref{fig:error_nlev}, we see that there are no problems with the accuracy of the numerical solution, regardless of the formulation used.

\begin{figure}[!tb]
  \centering
  \includegraphics[width=0.95\linewidth]{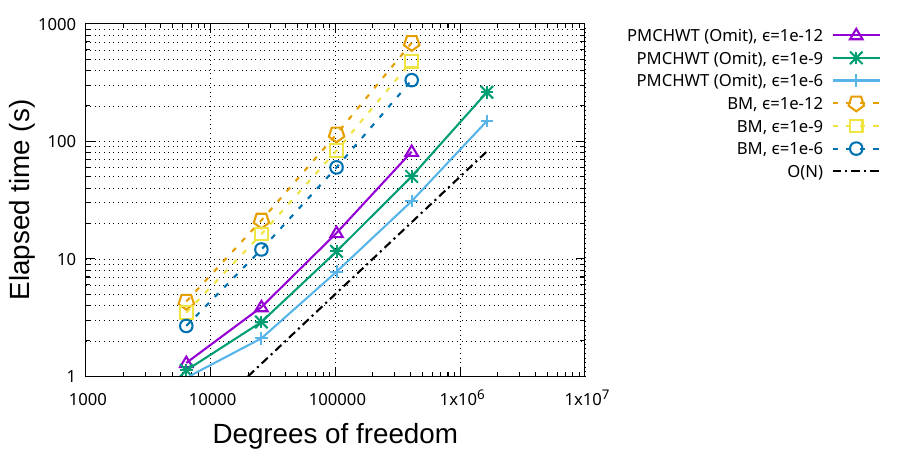}
  \caption{Computational time of the fast direct solver based on PMCHWT (Omit) and BM formulations.
      The labels ``PMCHWT (Omit)'' and ``BM'' correspond to the fast direct solver based on \eqref{eq:pmchwt_omit_bie} and \eqref{eq:bm_bie}, respectively.
      In the labels, $\epsilon$ represents the threshold of the column-pivoted QR decomposition in the fast direct solver.
      The results of the PMCHWT (Omit) formulation are represented by solid lines, while those of the BM formulation are shown as dashed lines.
      In this figure, the angular frequency $\omega$ and the material constants $\varepsilon^\pm$ are fixed.
      The PMCHWT (Omit) formulation is faster than BM formulation.
  }
  \label{fig:time_pmOmit_bm}
\end{figure}
\begin{figure}[!tb]
  \centering
  \includegraphics[width=0.95\linewidth]{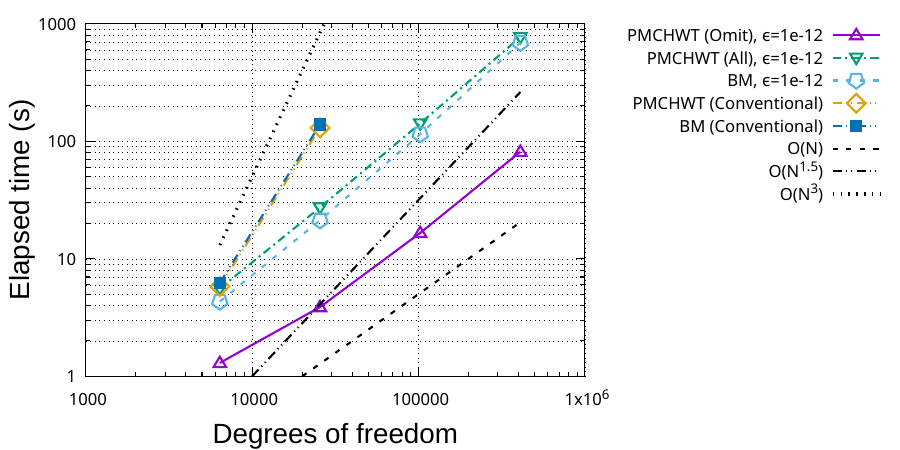}
  \caption{Computational time at threshold $\epsilon = 10^{-12}$. 
      The labels ``PMCHWT (Omit)'',  ``PMCHWT (All)'' and ``BM'' correspond to the fast direct solver based on \eqref{eq:pmchwt_omit_bie}, \eqref{eq:pmchwt_all_bie} and \eqref{eq:bm_bie}, respectively.
      In the labels, $\epsilon$ represents the threshold of the column-pivoted QR decomposition in the fast direct solver.
      The labels ``PMCHWT (Conventional)'',  ``BM (Conventional)'' correspond to \eqref{eq:pmchwt_omit_bie} and \eqref{eq:bm_bie}, respectively, and are solved by a standard partial pivoted LU decomposition.
      The results of the PMCHWT (Omit) formulation are represented by solid lines, while those of the other formulations are shown as dashed lines.
      The computational complexity of the proposed fast direct solvers is bounded by $O(N^{1.5})$.
  }
  \label{fig:time_eps12}
\end{figure}
\begin{figure}[!tb]
  \centering
  \includegraphics[width=0.95\linewidth]{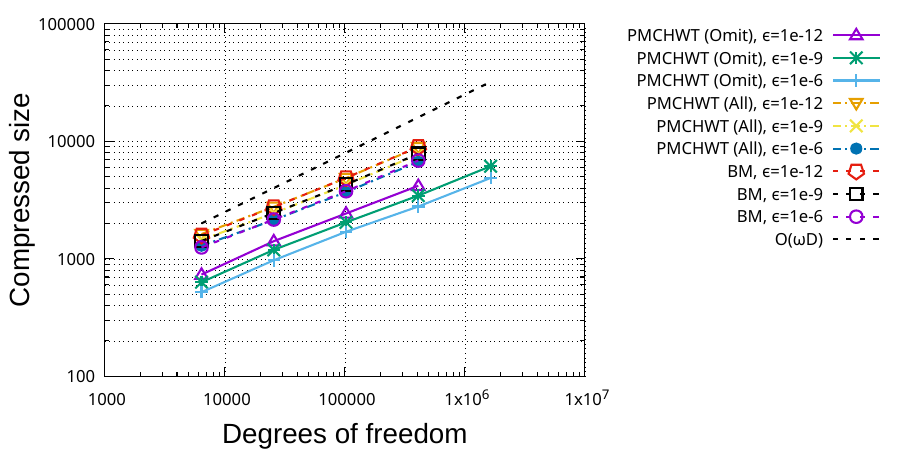}
  \caption{
    Degrees of freedom of original system versus compressed system with respect to increasing the number of inclusions.
      The labels ``PMCHWT (Omit)'',  ``PMCHWT (All)'' and ``BM'' correspond to the fast direct solver based on \eqref{eq:pmchwt_omit_bie}, \eqref{eq:pmchwt_all_bie} and \eqref{eq:bm_bie}, respectively.
      In the labels, $\epsilon$ represents the threshold of the column-pivoted QR decomposition in the fast direct solver.
      The results of the PMCHWT (Omit) formulation are represented by solid lines, while those of the other formulations are shown as dashed lines.
      The size of the compressed system roughly scales with $O(\omega D)$, where $\omega$ is a fixed angular frequency and $D$ denotes the diameter of the entire inclusion set.
  }
  \label{fig:size_nlev}
\end{figure}

\begin{table}[H]
  \centering
  \caption{Speedup ratio of the PMCHWT (Omit) formulation relative to the BM formulation.
      The former is about six times faster than the latter.
  }
  \resizebox{1.1\textwidth}{!}{%
    \begin{tabular}{rrrrcrrrcrrr}
    \toprule
    \multicolumn{1}{c}{DOF} & \multicolumn{3}{c}{$\epsilon = 10^{-12}$} & & \multicolumn{3}{c}{$\epsilon = 10^{-9}$} & & \multicolumn{3}{c}{$\epsilon = 10^{-6}$} \\
    \cmidrule{2-4} \cmidrule{6-8} \cmidrule{10-12}
    & \multicolumn{2}{c}{Time (s)}   & Speedup & & \multicolumn{2}{c}{Time (s)} & Speedup & & \multicolumn{2}{c}{Time (s)}   & Speedup \\
    \cmidrule{2-3} \cmidrule{6-7} \cmidrule{10-11}
    & BM & PMCHWT (Omit)  &  & & BM & PMCHWT (Omit) &  & & BM & PMCHWT (Omit) &  \\
    \midrule
    6400   & 4.33   & 1.30  & 3.33  & & 3.49   & 1.12  & 3.10  & & 2.68   & 0.97  & 2.76 \\
    25600  & 21.26  & 3.88  & 5.48  & & 16.22  & 2.90  & 5.59  & & 12.01  & 2.11  & 5.70 \\
    102400 & 114.98 & 16.59 & 6.93  & & 83.24  & 11.57 & 7.20  & & 59.97  & 7.77  & 7.72 \\
    409600 & 689.42 & 81.49 & 8.46  & & 479.48 & 50.22 & 9.55  & & 331.20 & 30.95 & 10.70 \\
    \midrule
    Average &        &       & 6.05  & &        &       & 6.36  & &        &       & 6.72 \\
    \bottomrule
    \end{tabular}%
  }
  \label{tb:speed}
\end{table}
On the other hand, Figures \ref{fig:time_pmOmit_bm} and \ref{fig:time_eps12} indicate that the fast direct solver based on the PMCHWT (Omit) formulation is more advantageous in terms of computational efficiency.
For the same threshold, Figure \ref{fig:time_pmOmit_bm} shows that the fast direct solver based on the PMCHWT (Omit) formulation is approximately six times faster than those based on the BM formulations as summarized in Table \ref{tb:speed}.
Although the speed difference is not very large, we can also see from Figure \ref{fig:time_eps12} that the fast direct solver based on the BM formulation outperforms that of the PMCHWT (All) formulation at the same threshold.
From the same figure, it can be observed that the computational complexity of the fast direct solver is bounded by $O(N^{1.5})$, which is consistent with our expectations.
Figure \ref{fig:size_nlev} shows that degrees of freedom of the compressed system is also bounded by $O(\omega D)$, where $D$ is the diameter of the bounded domain.
Moreover, it can be seen from this figure that the PMCHWT (Omit) formulation reduces the size of the compressed system to approximately half that of the BM formulation as summarized in Table \ref{tb:size_ratio}.
In the former formulation, only the exterior layer potentials appear in the off-diagonal blocks of the block linear system \eqref{eq:block_linear}. This likely accounts for the reduced number of skeletons required for low-rank approximation.
\begin{table}[H]
  \centering
  \caption{Compression ratio of the PMCHWT (Omit) formulation relative to the BM formulation.
      The PMCHWT based approach achieves a 50\% greater reduction in the compressed system's degrees of freedom compared
      to the solver based on the Burton--Miller formulation.
  }
  \resizebox{1.1\textwidth}{!}{%
    \begin{tabular}{rrrrcrrrcrrr}
    \toprule
    \multicolumn{1}{c}{DOF} & \multicolumn{3}{c}{$\epsilon = 10^{-12}$} & & \multicolumn{3}{c}{$\epsilon = 10^{-9}$} & & \multicolumn{3}{c}{$\epsilon = 10^{-6}$} \\
    \cmidrule{2-4} \cmidrule{6-8} \cmidrule{10-12}
    & \multicolumn{2}{c}{Compressed size}   & Ratio & & \multicolumn{2}{c}{Compressed size} & Ratio & & \multicolumn{2}{c}{Compressed size}   & Ratio \\
    \cmidrule{2-3} \cmidrule{6-7} \cmidrule{10-11}
    & BM & PMCHWT (Omit)  &  & & BM & PMCHWT (Omit) &  & & BM & PMCHWT (Omit) &  \\
    \midrule
    6400   & 1584 & 734  & 0.463  & & 1418 & 632  & 0.446  & & 1250 & 520  & 0.416 \\
    25600  & 2760 & 1412 & 0.511  & & 2452 & 1192 & 0.486  & & 2158 & 972  & 0.450 \\
    102400 & 4878 & 2430 & 0.498  & & 4284 & 2044 & 0.477  & & 3758 & 1698 & 0.452 \\
    409600 & 8972 & 4180 & 0.466  & & 7838 & 3448 & 0.440  & & 6858 & 2786 & 0.406 \\
    \midrule
    Average &        &       & 0.485  & &        &       & 0.462  & &        &       & 0.431 \\
    \bottomrule
    \end{tabular}%
  }
  \label{tb:size_ratio}
\end{table}

\subsection{Testing various frequencies} \label{sec:omega}
In the previous subsection, we investigated the scaling of the compressed system's size as the diameter of the bounded domain increases.
However, it is also necessary to examine how the compressed size behaves when the frequency is varied.
The calculation conditions used in this subsection are shown in Table \ref{tb:omega_condition}.
The number of inclusions is set to 64, as this is the maximum size that can be computed by the conventional BIEM based on the BM formulation.
\begin{table}[h]
    \centering
    \caption{Calculation conditions used in Section \ref{sec:omega}.}
    \label{tb:omega_condition}
    \begin{tabular}{lr}
    \toprule
      Parameter & Value \\
      \midrule
      Angular frequency ($\omega$)    & \{1, 2, 4, 8\} \\
      Material constant of exterior domain ($\varepsilon^{+}$) & 1 \\
      Material constant of interior domain ($\varepsilon^{-}$) & 5 \\
      Number of inclusions & 64 \\
      Number of quadrature points per inclusion & 200 \\
      Threshold of the column-pivoted QR decomposition ($\epsilon$) & \{$10^{-12}, 10^{-9}, 10^{-6}$\} \\
    \bottomrule
    \end{tabular}
\end{table}

\begin{figure}[!tb]
  \centering
  \includegraphics[width=0.95\linewidth]{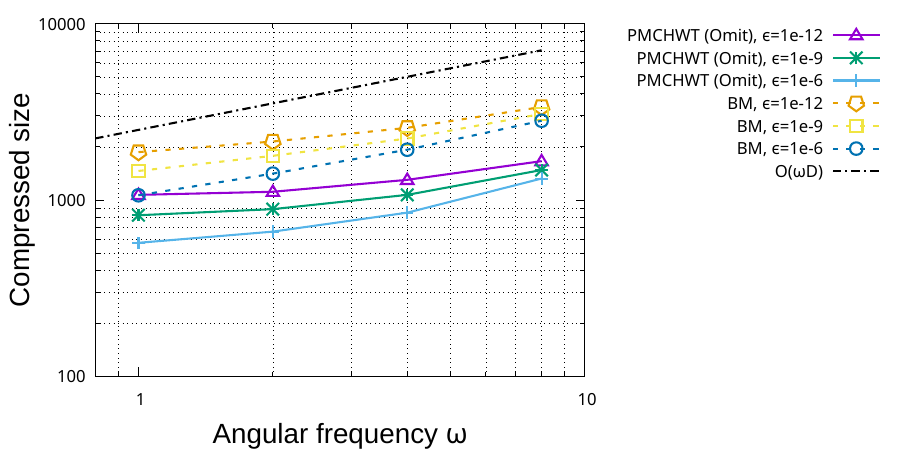}
  \caption{
    Degrees of freedom of compressed system with respect to increasing the frequency $\omega$ for a fixed degrees of freedom 25600.
      The labels ``PMCHWT (Omit)'' and ``BM'' correspond to the fast direct solver based on \eqref{eq:pmchwt_omit_bie} and \eqref{eq:bm_bie}, respectively.
      In the labels, $\epsilon$ represents the threshold of the column-pivoted QR decomposition in the fast direct solver.
      The results of the PMCHWT (Omit) formulation are represented by solid lines, while those of the BM formulation are shown as dashed lines.
      The size of the compressed system is well bounded by $O(\omega D)$, where $D$ denotes the diameter of the entire inclusion set.
  }
  \label{fig:size_omega}
\end{figure}
\begin{figure}[!tb]
  \centering
  \includegraphics[width=0.95\linewidth]{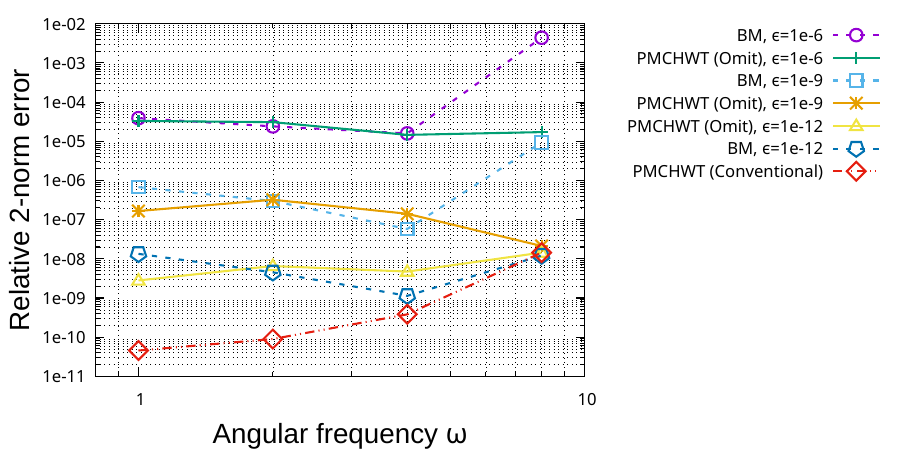}
  \caption{Relative 2-norm error for various frequencies for a fixed degrees of freedom 25600.
      The labels ``PMCHWT (Omit)'' and ``BM'' correspond to the fast direct solver based on \eqref{eq:pmchwt_omit_bie} and \eqref{eq:bm_bie}, respectively.
      In the labels, $\epsilon$ represents the threshold of the column-pivoted QR decomposition in the fast direct solver.
      The results of the PMCHWT (Omit) formulation are represented by solid lines, while those of the BM formulation are shown as dashed lines.
      The results of BM (conventional) is used as the reference.
      This figure indicates that the proposed fast direct solvers maintain accuracy within the low-frequency range.
  }
  \label{fig:error_omega}
\end{figure}
\begin{figure}[!tb]
  \centering
  \includegraphics[width=0.95\linewidth]{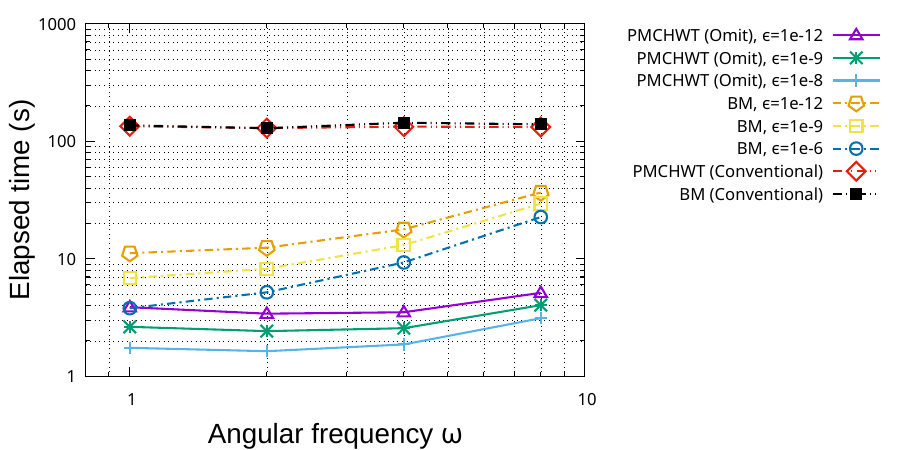}
  \caption{Computational time for various frequencies for a fixed degrees of freedom 25600.
      The labels ``PMCHWT (Omit)'' and ``BM'' correspond to the proposed fast direct solver based on \eqref{eq:pmchwt_omit_bie} and \eqref{eq:bm_bie}, respectively.
      In the labels, $\epsilon$ represents the threshold of the column-pivoted QR decomposition in the fast direct solver.
      The results of the PMCHWT (Omit) formulation are represented by solid lines, while those of the BM formulation are shown as dashed lines.
      The PMCHWT (Omit) based fast direct solver remains consistently more efficient than the BM based solver
  }
  \label{fig:time_omega}
\end{figure}

Figure \ref{fig:size_omega} illustrates the scaling of the compressed system's size with respect to the frequency $\omega$.
By combining the results from Figures \ref{fig:size_nlev} and \ref{fig:size_omega},
we conclude that the size of the system compressed by the proposed solver scales as $O(\omega D)$.
As in the previous subsection, we further investigate the accuracy and computational efficiency of the proposed method.
Figures \ref{fig:error_omega} and \ref{fig:time_omega} represents the relative 2-norm error and the computational time.
The solution obtained from the Burton--Miller formulation solved by a standard partial pivoted LU decomposition is employed as the reference of the relative 2-norm error.
The former figure indicates that numerical accuracy is generally maintained within the low-frequency range examined, with the exception of the BM-based fast direct solver at a frequency $\omega = 8$ with the thresholds $\epsilon = 10^{-6}$ and $10^{-9}$.
The BM formulated fast direct solver necessitates the calculation of all interior layer potentials, which are inherently redundant. This requirement not only increases the computational load but may also lead to a degradation in numerical accuracy.
The latter figure demonstrates that the PMCHWT (Omit) based fast direct solver remains consistently more efficient than the BM based solver, even as the frequency varies.


\subsection{Discussions}
  The numerical examples in Section \ref{sec:verify}--\ref{sec:omega} may depend strongly on a grid arrangement.
  We discuss the modifications required for applying the proposed solver to more general scatterer arrangement, and how these changes are expected to influence the results.
  In arbitrary configurations, it is necessary to implement a flexible tree-based clustering of integration points for scatterers, which is crucial for parallel computing efficiency.
  Currently, parallel efficiency is optimized for grid-aligned scatterers by indexing and clustering according to the Morton Z-order curve (Sagan, 2012).
  However, to accommodate arbitrary scatterer distributions, more flexible approaches such as k-means clustering for the integration points will be required (Ma et al., 2022).
  In this clustering approach, the distance between cell centers is employed to determine which cell pairs are assigned to near-field interactions within the proxy method.
  A cell pair is classified as near if the distance between them, scaled by a constant $\eta > 0$, is smaller than the diameter of either cell.
  This constant $\eta$ corresponds to the parameter for the admissibility condition (Bebendorf, 2000).
  Under such a clustering approach, scatterers themselves can no longer be treated as cells at the leaf level of a tree, and an appropriate data structure must be designed into the code.
  For example, higher-order discretization requires precise identification of the correspondence between quadrature points and their respective scatterers or cells, thereby increasing the overall implementation cost.

  Despite the challenges in implementation, the aforementioned approach enables the proposed solver to handle multiple scattering problems involving randomly distributed scatterers, including scatterers that are close to each other, provided a sufficiently fine hierarchical clustering is constructed.
  This capability is important for analyzing metamaterials and heterogeneous media.
  Regarding the computational efficiency in such scenarios, two primary issues arise.
  First, parallelization efficiency may be slightly compromised.
  To optimize the parallel performance, the number of quadrature points per cell after clustering should be uniform across each level of the tree structure.
  However, such uniformity is difficult to achieve with arbitrary distributions.
  Second, while $O(N^{1.5})$ scaling is expected to hold for scatterers that are distributed randomly but uniformly on average, performance may degrade if the distribution is highly non-uniform.
  In such cases, the number of levels available for multi-level compression may decrease, leading to reduced computational efficiency.
  As a special case, it is known that $O(N)$ scaling can be achieved with if a scatterer has elongated structures (Martinsson and Rokhlin, 2007).
  Following this logic, the proposed solver may achieve $O(N)$ complexity when scatterers are distributed in a thin, strip-like arrangement.

\section{Concluding remarks}
In this paper, we developed a fast direct solver using boundary integral equations for two-dimensional Helmholtz transmission problems involving multiple inclusions.
The proposed solver provides a promising alternative to iterative methods, as the latter often suffer from slow convergence in the presence of numerous scatterers.
One of the findings in this work is that the PMCHWT formulation can be more effective than the Burton--Miller approach when integrated with a proxy-based fast direct solver.
By omitting the interior integral representation, the PMCHWT based solver achieved a six-fold speedup compared to the Burton--Miller formulation. Furthermore, the PMCHWT formulation reduced the size of the compressed system by half relative to its Burton--Miller counterpart under identical computational settings. These findings highlight the importance of selecting an appropriate boundary integral equation to maximize the performance of fast direct solvers in multiple scattering scenarios.
Moreover, the numerical results demonstrate that the proposed solver effectively compresses the system of linear algebraic equations to a size of $O(\omega D)$. Here, $\omega$ is the frequency  of the incident wave and $D$ is the diameter of the (smallest) bounding box enclosing the multiple inclusions.
Specifically, when the inclusions are arranged on a grid, the total computational cost scales at most as $O(N^{1.5})$ for a fixed frequency, confirming the efficiency of the approach for large-scale problems.

The proposed approach has certain limitations that suggest directions for future research. In this study, we maximized the efficiency of the PMCHWT formulation by identifying each individual inclusion as a leaf-level cell in the fast solver. While effective, this strategy is primarily limited to two-dimensional problems at relatively low frequencies.
For high-frequency or three-dimensional problems, a single inclusion cannot be accurately discretized with a small number of degrees of freedom, even when employing high-order discretization methods. To address this, it is necessary to introduce a hierarchical tree structure within each scatterer.
Such a scheme would require a switching strategy for the proxy-based interaction calculations: the solver initially account for both internal and external contributions within an inclusion, subsequently transitioning to use of only external layer potentials for disjoint inclusions.
Furthermore, the low-rank approximation condition should be updated from weak admissibility to strong admissibility to maintain efficiency and accuracy in these more complex scenarios.

\FloatBarrier

\section*{Acknowledgments}
This work was supported by
Japan Society for the Promotion of Science under KAKENHI grant number 24K20783,
and by the computational resources
of Camphor3 at the Kyoto University
provided through the projects 
``Joint Usage/Research Center for Interdisciplinary Large-scale Information Infrastructures (JHPCN)''
and ``High Performance Computing Infrastructure (HPCI)'' in Japan (project IDs: jh250045 and jh250070).

\end{document}